\documentclass[reqno]{amsart}

\usepackage{article-preamble}

\hypersetup{colorlinks=true,
    linkcolor=blue,
    citecolor=blue,
    pdftitle={A singular integral identity for surface measure},
    pdfauthor={Ryan E. G. Bushling}}

\definecolor{verylightgray}{rgb}{0.875,0.875,0.875}
\newcommand\normvec[4]{\draw[very thick,->] (#1,#2) -- ({#1+#3*cos(#4)},{#2+#3*sin(#4)});}

\begin{document}

\title[]{A singular integral identity for surface measure}
\author[]{Ryan E. G. Bushling}
\address{Department of Mathematics \\ University of Washington, Box 354350 \\ Seattle, WA 98195-4350}
\email{reb28@uw.edu}

\maketitle

\vs{-1}

\begin{abstract}
    We prove that the integral of a certain Riesz-type kernel over $(n-1)$-rectifiable sets in $\R^n$ is constant, from which a formula for surface measure immediately follows. Geometric interpretations are given, and the solution to a geometric variational problem characterizing convex domains follows as a corollary, strengthening a recent inequality of Steinerberger.
\end{abstract}

\section{Introduction and main results}

In \cite{steinerb2022inequality}, Steinerberger proves an inequality inspired by the following simple observation: if $\Omega$ is a smoothly bounded convex domain and $x,y \in \partial\Omega$ are close with respect to the Euclidean distance, then the normal vectors at $x$ and $y$ are nearly orthogonal to $x-y$, where the measure of ``closeness" hinges on the curvature of $\partial\Omega$. Leveraging this from a probabilistic standpoint, he concludes the following. Let $\mathcal{H}^{n-1}$ denote $(n-1)$-dimensional Hausdorff measure and $\alpha_{n-1} := \mathcal{L}^{n-1}(B(0,1))$ the Lebesgue measure of the unit ball in $\R^{n-1}$. \vs{-0.15}

\begin{prop} \label{prop:steinerb}
   For every bounded, $C^1$-bounded domain $\Omega \subset \R^n$ with outward unit normal vector field $\nu$,
    \begin{equation*}
        \int_{\partial\Omega} \int_{\partial\Omega} \frac{|\langle x-y, \nu(y) \rangle \langle x-y, \nu(x) \rangle|}{\| x-y \|^{n+1}} \, d\mathcal{H}^{n-1}(y) \+ d\mathcal{H}^{n-1}(x) \geq \alpha_{n-1} \mathcal{H}^{n-1}(\partial\Omega).
    \end{equation*}
    Moreover, equality holds if and only if $\Omega$ is convex.
\end{prop}

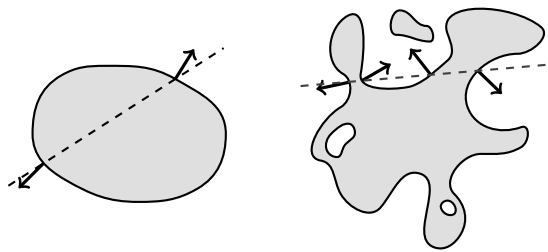
\begin{figure}[h!]
    \begin{tikzpicture}[scale=0.8]

    \begin{scope}[xshift=-40,yshift=60,scale=0.73]
    
    \draw[thick,fill=verylightgray] (0.15,-1.55) to [out=-180,in=-30] (-1.5,-1.05) to [out=150,in=-90] (-2.25,0) to [out=90,in=-150] (-1.6,1.33) to [out=30,in=180] (-0.6,1.55) to [out=0,in=150] (1,1.25) to [out=-30,in=90] (2.15,0) to [out=-90,in=30] (1.5,-1.25) to [out=-150,in=0] (0.15,-1.55);
    
    \normvec{-2}{-0.67}{0.8}{-135}
    \normvec{1}{1.25}{0.8}{58}
    \draw[thick,dashed] (-2.8,-1.182) -- (2.1,1.954);
    
    \end{scope}
    
    \begin{scope}[xshift=40,scale=0.37]
    
    \filldraw[thick,fill=verylightgray](3.6,2.4) 
    to[out=50,in=160] (5.4,3.9) 
    to[out=-20,in=80] (5.7,2.4) 
    to[out=-100,in=160] (6,.6)  
    to[out=-20,in=-120] (7.2,1.2) 
    to[out=60,in=-45]  (6.9,3) 
    to[out=135,in=180]  (8.1,4.8) 
    to[out=0,in=-150] (9.9,5.1) 
    to[out=30,in=20, looseness=1.3]  (9.9,6) 
    to[out=-160,in=-70] (7.5,6.9) 
    to[out=110,in=180]  (9.75,9.3) 
    to[out=0,in=-40] (10.8,10.5) 
    to[out=140,in=110]  (6.6,10.2) 
    to[out=-70,in=40] (6.6,9) 
    to[out=-140,in=20,looseness=0.9] (5,7.85)
    to[out=-160,in=10,looseness=0.9] (4.8,7.8) 
    to[out=-170,in=-100,looseness=0.9] (2.7,8.4) 
    to[out=80,in=30] (2.1,11) 
    to[out=-150,in=150] (1.9,8.4) 
    to[out=-30,in=30] (1.8,6.9) 
    to[out=-150,in=170] (0.9,4.5) 
    to[out=-10,in=120] (1.8,3) 
    to[out=-60,in=-130] (3.6,2.4);
    
    \filldraw[thick,fill=verylightgray,xshift=40,yshift=10] (2.9,9.8) 
    to[out=0,in=180] (4.4,9.5)
    to[out=0,in=0] (4.25,10.2)
    to[out=180,in=0] (3.05,10.9)
    to[out=180,in=180] (2.9,9.8);
    
    \filldraw[thick,fill=white] (6.3,2.4) 
    to[out=-80,in=-145] (6.9,2.1)
    to[out=45,in=0] (6.6,2.7)
    to[out=180,in=100](6.3,2.4);
    
    \filldraw[thick,fill=white,rotate=0,xshift=-33,yshift=50] (2.3,3.4) 
    to[out=90,in=-140] (2.7,4.1)
    to[out=40,in=180] (3.3,4.4)
    to[out=0,in=90] (3.6,4)
    to[out=-90,in=45] (3.2,3.5)
    to[out=-135,in=70] (3,3.1)
    to[out=-110,in=0] (2.8,2.9)
    to[out=180,in=-90] (2.3,3.4);
    
    \normvec{2.175}{8.02}{1.5}{-168}
    \normvec{2.75}{8.09}{1.5}{30}
    \normvec{5.85}{8.35}{1.5}{128}
    \normvec{7.96}{8.54}{1.5}{-44}
    
    \draw[thick,dashed,darkgray] (0,7.85) -- (11,8.8);
    
    \end{scope}
    
    \end{tikzpicture}
    \caption{For a $C^1$-bounded convex region $\Omega$, the line segment between any two points $x,y \in \partial\Omega$ is such that $\langle x-y, \nu(y) \rangle \langle x-y, \nu(x) \rangle \leq 0$, and this quantity vanishes quickly as $\| x-y \| \to 0$. If $\Omega$ is not convex, then $\langle x-y, \nu(y) \rangle \langle x-y, \nu(x) \rangle$ can be positive and there may be nearby points $x,y \in \partial\Omega$ such that $\langle x-y, \nu(y) \rangle \langle x-y, \nu(x) \rangle$ is not small relative to $\| x-y \|$.}
    \label{fig:domains}
\end{figure}

\begin{figure}[b]
\centering
    \begin{tikzpicture}[scale=0.7]
    
    \draw[thick,dashed,darkgray] (2.25,-3.6328) -- (-3.6,2.665);
    
    \draw[very thick, looseness=1.2] (0,-0.6) 
        to[out=180, in=100] (-0.8,-1.5) 
        to[out=-80, in=0] (-2,-2.8) 
        to[out=180, in=-90] (-3.5,-1) 
        to[out=90, in=-135] (-2.75,0.1) 
        to[out=45, in=-105] (-2.4,0.7) 
        to[out=75, in=155] (-1.2,1.1) 
        to[out=-25, in=-150] (0.4,1.4) 
        to[out=30, in=90] (2,1) 
        to[out=-90, in=0] (1.45,0.25) 
        to[out=180, in=180] (1.45,0.8) 
        to[out=0, in=90] (2.2,0) 
        to[out=-90, in=0] (1.2,-1.1) 
        to[out=180, in=-30] (0.6,-.85) 
        to[out=150, in=0] (0,-0.6);
    \draw[very thick] (-2,-1.4)
        to[out=-130, in=0] (-2.35,-1.6)
        to[out=180, in=-100] (-2.5,-1.25)
        to[out=80, in=-100] (-1.6, -.5)
        to[out=80, in=0] (-1.7, -.3)
        to[out=180, in=90] (-2, -.55)
        to[out=-90, in=50] (-2,-1.4);
    
    \draw[very thick](1.6,-2) 
    to[out=-30, in=80] (1.8,-2.6) 
    to[out=-100, in=-20] (1.2, -2.9) 
    to[out=160, in=-110] (1.1, -2.4)
    to[out=70, in=150] (1.6,-2);
    
    \draw[very thick](-2.9,.3) 
        to[out=160, in=-140](-2.8,1.4)
        to[out=40, in=-60](-2.6, 2)
        to[out=120, in=180](-2, 3.1)
        to[out=0, in=90](-1.4, 2.5)
        to[out=-90, in=10](-1.7, 2.1)
        to[out=-180, in=90](-2, 1.5) 
        to[out=-90, in=60](-2.1, .7) 
        to[out=-120, in=-10](-2.9,.3);
    
    \draw[very thick](-1.7,1.8) 
        to[out=135, in=-120] (-1.9,2.6) 
        to[out=60, in=100] (-1.3, 2)
        to[out=-80, in=-45] (-1.7,1.8);
    
    \draw[very thick] (2.4,2)
        to[out=-130, in=-130, looseness=2] (2.9,1.6)
        to[out=50, in=-45] (3.3, 3)
        to[out=135, in=80] (2.8, 2.9)
        to[out=-100, in=140] (3.4, 1.7)
        to[out=-40, in=-50] (3.8, 2.3)
        to[out=130, in=50] (2.4,2);
    
    \normvec{1.56}{-2.89}{0.68}{-66}
    \normvec{1.1}{-2.4}{0.68}{157}
    \normvec{-0.47}{-0.7}{0.68}{-61}
    \normvec{-2.01}{0.95}{0.68}{-14}
    \normvec{-2.15}{1.1}{0.68}{125}
    \normvec{-2.62}{1.6}{0.68}{147}
    
    \end{tikzpicture}
    \caption{The normal vectors to this immersed submanifold are oriented in such a way that the signs of the angles formed with any line sum to $0$. This motivates the definition of the ``orientation cancellation condition."}
    \label{fig:cells}
\end{figure}
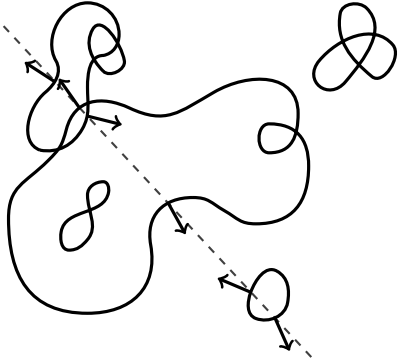

What prevents the inequality from being an equality in general is the absolute value: for $\Omega$ open and $C^1$-bounded, the sign of $\langle x-y, \nu(y) \rangle \langle x-y, \nu(x) \rangle$ is constant precisely when $\Omega$ is convex, and dropping the absolute value results in a ``systematic cancellation" that turns the inequality into a formula for surface measure (cf.~Figure \ref{fig:domains}). This remedy begs the question whether boundaries of domains are the natural class of hypersurface with which to work in this context, as the setup only requires a normal vector field that is distributed ``consistently" across the surface, as in Figure \ref{fig:cells}. In \S \ref{ss:occ}, we specify such a class of surface/vector field pairs and say that its members satisfy the \textit{orientation cancellation condition}. The class includes all boundaries of bounded, $C^1$-bounded domains and all compact, oriented, immersed smooth $(n-1)$-manifolds (both with their outward unit normal vector fields), as well as a host of lower-regularity sets with vector fields that do not arise as the result of an ``orientation."

In this setting, we can prove the following theorem. For the remainder of this section and subsequently in \S \ref{s:proof}, $\Sigma \subset \R^n$ denotes an $(n-1)$-rectifiable set and $\nu$ a measurable unit normal vector field on $\Sigma$ (cf.~\S\ref{ss:gmt}).

\begin{thm} \label{thm:main}
    For every $(\Sigma,\nu)$ satisfying the orientation cancellation condition (cf.~\S \ref{ss:occ}), the identity
    \begin{equation} \label{eq:main-1}
        \int_\Sigma \frac{\langle x-y, \nu(y) \rangle \langle y-x, \nu(x) \rangle}{\| x-y \|^{n+1}} \+ d\mathcal{H}^{n-1}(y) = \alpha_{n-1}.
    \end{equation}
    holds for $\mathcal{H}^{n-1}$-a.e.~$x \in \Sigma$. Consequently,
    \begin{equation} \label{eq:main-2}
        \frac{1}{\alpha_{n-1}} \int_\Sigma \int_\Sigma \frac{\langle x-y, \nu(y) \rangle \langle y-x, \nu(x) \rangle}{\| x-y \|^{n+1}} \, d\mathcal{H}^{n-1}(y) \+ d\mathcal{H}^{n-1}(x) = \mathcal{H}^{n-1}(\Sigma). \vs{0.15}
    \end{equation}
\end{thm}

\begin{figure}[t]
\vs{0.15}
\centering
    \begin{tikzpicture}[scale=0.6]
    
    \draw[fill=verylightgray,smooth,looseness=0.9] (0,-2) to [out=-80,in=55] (-0.2,-3.4) to [out=-125,in=-15] (-1.75,-3.6) to [out=165,in=-40] (-4,-2.3) to [out=140,in=-90] (-5.6,0) to [out=90,in=-135] (-5,1.3) to [out=45,in=180] (-4,1.7) to [out=0,in=170] (-2.5,1.4) to [out=-10,in=-90] (0.15,2.5) to [out=90,in=-55] (-0.5,3.6) to [out=125,in=-90] (-1.3,5) to [out=90,in=180] (0,5.8) to [out=0,in=150] (2.3,5) to [out=-30,in=90] (3,4) to [out=-90,in=60] (2.3,3) to [out=-120,in=135] (2.5,1.6) to [out=-45,in=180] (4,1.25) to [out=0,in=160] (5,1.1) to [out=-20,in=90] (5.6,0.4) to [out=-90,in=40] (5,-0.8) to [out=-140,in=0] (3.6,-1.35) to [out=180,in=-35] (2.6,-0.88) to [out=145,in=0] (1,0) to [out=180,in=45] (0.2,-0.45) to [out=-135,in=90] (-0.2,-1.2) to [out=-90,in=100] (0,-2);
    
    \draw[thick,smooth,looseness=0.9] (0,-2) to [out=-80,in=55] (-0.2,-3.4) to [out=-125,in=-10] (-1.15,-3.73);
    
    \draw[thick,smooth,looseness=0.9] (-2.35,-3.4) to [out=157,in=-40] (-4,-2.3) to [out=140,in=-90] (-5.6,0) to [out=90,in=-135] (-5,1.3) to [out=45,in=180] (-4,1.7) to [out=0,in=170] (-2.5,1.4) to [out=-10,in=-90] (0.15,2.5) to [out=90,in=-55] (-0.5,3.6) to [out=125,in=-90] (-1.3,5) to [out=90,in=180] (0,5.8) to [out=0,in=150] (2.3,5) to [out=-30,in=90] (3,4) to [out=-90,in=60] (2.3,3) to [out=-120,in=130] (2.43,1.68);
    
    \draw[thick,smooth,looseness=0.9] (5.13,1.04) to [out=-25,in=90] (5.6,0.4) to [out=-90,in=40] (5,-0.8) to [out=-140,in=0] (3.6,-1.35) to [out=180,in=-35] (2.6,-0.88) to [out=145,in=-40] (2.28,-0.63);
    
    \draw[thick,smooth,looseness=0.9] (1.3,-0.035) to [out=165,in=0] (1,0) to [out=180,in=45] (0.2,-0.45) to [out=-135,in=90] (-0.2,-1.2) to [out=-90,in=100] (0,-2);
        
    \foreach \ang/\len in {0/6.2,1/7,2/8,3/8.5,4/9,5/9,6/7.5,7/6.5,8/6,9/5,10/4.5,11/4.5,12/4.75,13/5.5} {
        \draw[thick,gray,dotted] (0,-2) -- ({\len*cos(5+\ang*180/7)}, {-2+\len*sin(5+\ang*180/7)});
    }
    
    \fill[pattern=north west lines,pattern color=lightgray,line width=1in,opacity=0.9] (0,-2) -- (6.018,1.575) to [out=100,in=-50] (4.424,4.666) -- cycle;
    \fill[pattern=north west lines,pattern color=lightgray,opacity=0.9] (0,-2) -- (-5.158,-5.064) to [out=-35,in=170] (-2.765,-6.166) -- cycle;
    
    \normvec{1.81}{-0.27}{1.1}{-127}
    \normvec{3.427}{1.26}{1.1}{90}
    \normvec{-1.688}{-3.606}{1.1}{-106}
    
    \draw[ultra thick,smooth,looseness=0.9] (-1.15,-3.73) to [out=170,in=-23] (-2.35,-3.4);
    \draw[ultra thick,smooth,looseness=0.9] (2.43,1.68) to [out=-50,in=135] (2.5,1.6) to [out=-45,in=180] (4,1.25) to [out=0,in=160] (5,1.1) to [out=-20,in=155] (5.13,1.04);
    \draw[ultra thick,smooth,looseness=0.9] (2.28,-0.63) to [out=140,in=-15] (1.3,-0.035);
    
    \draw[thick,darkgray] (0,-2) -- ({7*cos(5+1*180/7)}, {-2+7*sin(5+1*180/7)});
    \draw[thick,darkgray] (0,-2) -- ({-6*cos(5+1*180/7)}, {-2-6*sin(5+1*180/7)});
    
    \draw[thick,darkgray] (0,-2) -- ({8*cos(5+2*180/7)}, {-2+8*sin(5+2*180/7)});
    \draw[thick,darkgray] (0,-2) -- ({-5*cos(5+2*180/7)}, {-2-5*sin(5+2*180/7)});

    \filldraw (0,-2) circle (0.075cm);
    \filldraw[white] (0.125,-2.55) rectangle (0.625,-2.05);
    \node[below right] at (0,-2) {$x$};
    
    \end{tikzpicture}
    \vs{0.15}
    \caption{The proof of the Theorem formalizes the following idea: if $\Sigma$ (depicted here as the boundary of a $C^1$-bounded region) is partitioned into double cones with vertex at $x$, then each piece of $\Sigma$ that slices the double cone contributes approximately the same mass to the integral in Equation \eqref{eq:main-1}, up to a sign. (The weight factor $\langle x-y, \nu(y) \rangle/\| x-y \|^n$ is chosen precisely to guarantee this.) The singleton $\{ x \}$ has no mass, so the OCC implies that the contribution to the integral from this double cone is approximately the area of a single slice of the cone that is unit distance from $x$ and orthogonal to the axis of the cone.}
    \label{fig:polar}
\end{figure}
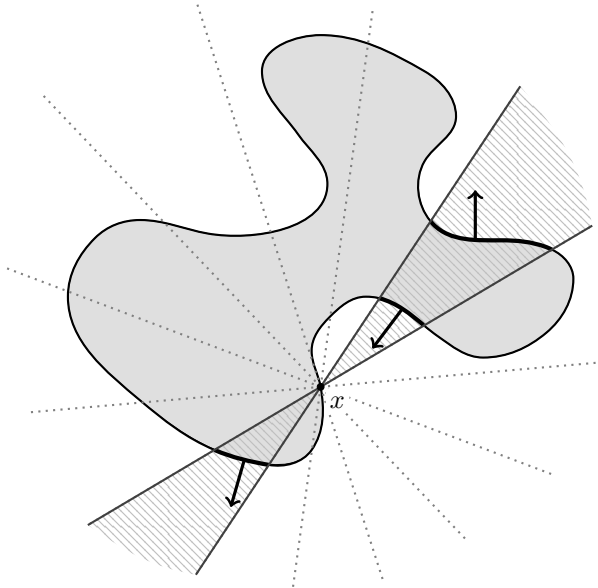

In plain language, the Theorem states the following: if $\Sigma$ were semitransparent, then the amount of $\Sigma$ that one would see while standing on the surface\textemdash counting each piece of $\Sigma$ positively or negatively according to its orientation relative to the viewer\textemdash would not depend on the point at which one stood. In fact, this quantity does not even depend on the surface $\Sigma$: it is a universal constant depending only on the dimension $n$, and taking $\Sigma = \bbs^{n-1}$ gives the constant explicitly. It follows immediately that the surface area of $\Sigma$ is proportional to the integral over all $x \in \Sigma$ of the signed surface area one sees from the vantage point $x$. While this interpretation is not apparent from the theorem statement, the heuristic is salient in the proof. See also Figure \ref{fig:polar}.

A more concrete consequence of the proof is that Equation \eqref{eq:main-1} holds for every $x \in \Sigma$ at which the orientation cancellation condition is satisfied. In particular, if $\Sigma = \partial E \in C^1$ for some bounded open set $E$ and if $\nu$ is outward-pointing, then the equation holds for all $x \in \partial E$. However, even if $E$ is a bounded set of finite perimeter with Gauss-Green measure $\mu_E = \nu_E \+ \mathcal{H}^{n-1} \restrict \partial^*\-E$, the orientation cancellation condition is still satisfied with $(\Sigma,\nu) = (\partial^*\-E,\nu_E)$ at $\mathcal{H}^{n-1}$-a.e.~$x \in \partial^*\-E$. (See \S \ref{ss:gmt}.)

In view of this discussion (formalized in Lemma \ref{lem:finite-perimeter} below), the Theorem implies Steinerberger's proposition under a milder regularity hypothesis.

\begin{cor} \label{cor:convex-inequality}
    For every bounded set $E \subset \R^n$ of finite perimeter,
    \begin{equation} \label{eq:convex-inequality}
        \frac{1}{\alpha_{n-1}} \int_{\partial^*\-E} \int_{\partial^*\-E} \frac{|\langle x-y, \nu_E(y) \rangle \langle y-x, \nu_E(x) \rangle|}{\| x-y \|^{n+1}} \, d\mathcal{H}^{n-1}(y) \+ d\mathcal{H}^{n-1}(x) \geq \mathcal{H}^{n-1}(\partial^*\-E).
    \end{equation}
    Furthermore, there is equality if and only if $E$ is $\mathcal{L}^n$-equivalent to a convex set.
\end{cor}

Notice that the inner integral
\begin{equation*}
    \int_{\partial^*\-E} \frac{|\langle x-y, \nu(y) \rangle \langle y-x, \nu(x) \rangle|}{\| x-y \|^{n+1}} \, d\mathcal{H}^{n-1}(y)
\end{equation*}
is unstable under $L^\infty$ perturbations of $\partial^*\-E$, although it \textit{is} stable under $C^1$ perturbations. As such, the magnitude of this ``energy" relative to the measure of the reduced boundary provides an interesting metric for how ``close" a set is to being convex. Steinerberger \cite{steinerb2022quadratic} substantiates this idea with an application to a geometric variational problem. Put slightly differently, the Corollary may be interpreted as stating (to paraphrase \cite{steinerb2022inequality}) that the solution set to a certain nonlocal isoperimetric problem, considered over the family of bounded finite-perimeter sets, is the family of all convex domains of a given perimeter.

\section{Definitions} \label{s:defs}

This short section describes the objects of study. A full account of this background can be found in \cite{krantz2008geometric} and \cite{maggi2012sets}. The generality is not so great as to defy classical methods, yet is sufficient to include the boundaries of all convex domains and the variety of hypersurfaces suggested by Figure \ref{fig:cells}.

\subsection{Rectifiable sets} \label{ss:gmt}

A set $\Sigma \subseteq \R^n$ is \define{$\bm{k}$-rectifiable} if $\mathcal{H}^k(\Sigma) < \infty$ and there exist a countable family $\{ F_i \}_{i \in I}$ of Lipschitz maps $F_i \!: A_i \subseteq \R^k \to \R^n$ and an $\mathcal{H}^k$-null set $\Sigma_0 \subseteq \Sigma$ such that
\begin{equation*}
    \Sigma = \Sigma_0 \cup \bigcup_{i \in I} F_i(A_i).
\end{equation*}
If $\Sigma$ is $k$-rectifiable, then, at $\mathcal{H}^k$-a.e.~$y \in \Sigma$, there exists a unique \textit{approximate tangent space} $T_y \Sigma \subseteq \R^n$, which coincides with the classical tangent space when $\Sigma$ is a smooth hypersurface. We call $\nu \!: \Sigma \to \bbs^{n-1}$ a \define{measurable unit normal vector field} on $\Sigma$ if it is $\mathcal{H}^k$-measurable and $\nu(y)$ is orthogonal to $T_y \Sigma$ for $\mathcal{H}^k$-a.e.~$y \in \Sigma$ at which the approximate tangent space is uniquely defined. In this case, the real-valued function $y \mapsto \langle T(y), \nu(y) \rangle$ is a measurable function on $\Sigma$ whenever $T \!: \Sigma \to \R^n$ is measurable, and the integral
\begin{equation*}
    \int_\Sigma \langle T(y), \nu(y) \rangle \, d\mathcal{H}^k(y)
\end{equation*}
is well-defined when the function is also integrable. (In particular, the proof of the Theorem implies that the integrand in the theorem statement is integrable for a.e.~$x \in \Sigma$.) The main features of rectifiable sets that we use are the notion of a normal vector field and the applicability of the coarea formula (cf.~\cite{krantz2008geometric}).

The language of the theory of sets of finite perimeter comes to bear in the Corollary. If $E \subseteq \R^n$ is such a set, then there is an $(n-1)$-rectifiable set $\partial^*\-E \subseteq \partial E$, the \textit{reduced boundary} of $E$, on which the \textit{measure-theoretic outward unit normal vector field} $\nu_E \!: \partial^*\-E \to \R^n$ is defined. The set $E$ comes with a natural vector-valued measure $\mu_E$, the \textit{Gauss-Green measure} of $E$, that takes the form $\mu_E = \nu_E \+ |\mu_E| = \nu_E \+ \mathcal{H}^{n-1} \restrict \partial^*\-E$. All bounded convex sets are sets of finite perimeter. While \cite{maggi2012sets} is our primary reference, \cite{ambrosio2000functions} contains some more nuanced results that we shall need as well. \vs{-0.15}

\subsection{The orientation cancellation condition} \label{ss:occ}

For each $x \in \R^n$ and $\omega \in \bbs^{n-1}$, let $L_{x,\omega} := x + \Span \omega$ denote the line through $x$ with direction vector $\pm\omega$.

\textbf{Definition.} Given an $(n-1)$-rectifiable set $\Sigma \subset \R^n$ and a measurable unit normal vector field $\nu \!: \Sigma \to \bbs^{n-1}$, we say the \define{orientation cancellation condition} (or \define{OCC}) \define{is satisfied} at a point $x \in \Sigma$ if the equation
\begin{equation} \label{eq:occ}
    \sum_{y \+ \in \+ \Sigma \+ \cap \+ L_{x,\omega}} \sgn \, \langle \omega, \nu(y) \rangle = 0
\end{equation}
holds for $\mathcal{H}^{n-1}$-a.e.~$\omega \in \bbs^{n-1}$, where $\sgn(\,\cdot\,)$ is the signum function with $\sgn \, 0 := 0$. If the OCC is satisfied at $\mathcal{H}^{n-1}$-a.e.~$x \in \Sigma$, we say that $(\Sigma,\nu)$ \define{satisfies the orientation cancellation condition}.

One can take this definition as an adaptation to lower-regularity sets of a concept from algebraic topology. An immersed smooth hypersurface $\Sigma \subset \R^n$ admitting a continuous normal vector field $\nu$ such that $(\Sigma,\nu)$ satisfies the OCC is said to have \textit{first Stiefel-Whitney class $0$}, and it is a theorem that this is equivalent to orientability. Naturally, the most salient example of such a surface is the boundary of a smoothly bounded open set, in which case $\sgn \, \langle \omega, \nu(y) \rangle$ alternates sign along successive values of $y \in \Sigma \cap L_{x,\omega}$. However, there are more general cell complexes with first Stiefel-Whitney class $0$ that are not oriented manifolds (cf.~\cite{milnor1974characteristic}). For example, there are many choices of orientation for the $1$-cells in Figure \ref{fig:cells} that give it first Stiefel-Whitney class $0$, although the resulting space need not admit a topology making it into an immersed, oriented submanifold of Euclidean space.

\section{Proofs of results} \label{s:proof}

We single out one computation before delving into the proof of the main theorem. \vs{-0.1}

\begin{lem} \label{lem:jacobian}
    Let $\Sigma \subset \R^n$ be an $(n-1)$-rectifiable set, $\nu \!: \Sigma \to \bbs^{n-1}$ a measurable unit normal vector field, $x \in \Sigma$ a point, and $\pi_x \!: \R^n \setminus \{ x \} \to \bbs^{n-1} \cong \partial B(x,1)$ the radial projection onto the unit sphere centered at $x$:
    \begin{equation*}
        \pi_x(y) := \frac{y-x}{\|y-x\|}.
    \end{equation*}
    Then the a.e.-defined Jacobian determinant $|J\pi_x| \!: \Sigma \setminus \{ x \} \to \R$ is given by
    \begin{equation*}
        |J\pi_x(y)| = \frac{|\langle x-y, \nu(y) \rangle|}{\| x-y \|^n}.
    \end{equation*}
\end{lem} \vs{-0.15}

\textit{Proof.} We employ the tensor notation of \cite{maggi2012sets}. Let $y \in \Sigma$ be a point at which $\nu$ is defined and let $(\mathbf{u}_i)_{i=1}^n$ be an orthonormal basis for $\R^n$ such that $\Span \, (\mathbf{u}_i)_{i=1}^{n-1} = T_y \Sigma$ and $\mathbf{u}_n = \nu(y)$. A routine computation gives the representation of the derivative $D_y \pi_x \!: \R^n \to \R^n$ in these coordinates:
\begin{equation*}
    D_y \pi_x = \sum_{i=1}^n \sum_{j=1}^n \frac{\delta_{ij} \| y-x \|^2 - (y_i - x_i)(y_j - x_j)}{\| y-x \|^3} \mathbf{u}_i \otimes \mathbf{u}_j,
\end{equation*}
where we write $z = \sum_{i=1}^n z_i \mathbf{u}_i$. The derivative at $y$ of the inclusion $\iota \!: \Sigma \hookrightarrow \R^n$ is \vs{-0.1}
\begin{equation*}
    D_y \iota = \sum_{j=1}^{n-1} \mathbf{u}_j \otimes \mathbf{u}_j, \vs{-0.1}
\end{equation*}
and composing with $D_y \pi_x$ gives the restriction of $D_y \pi_x$ to $T_y \Sigma$:
\begin{align*}
    D_y \pi_x |_{T_y \Sigma} &= D_y (\pi_x \circ \iota) = D_y \pi_x \circ D_y \iota \\
    &= \left( \sum_{i=1}^n \sum_{k=1}^n \frac{\delta_{ik} \| y-x \|^2 - (y_i - x_i)(y_k - x_k)}{\| y-x \|^3} \mathbf{u}_i \otimes \mathbf{u}_k \right) \! \left( \sum_{j=1}^{n-1} \mathbf{u}_j \otimes \mathbf{u}_j \right) \\
    &= \sum_{i=1}^n \sum_{j=1}^{n-1} \frac{\delta_{ij} \| y-x \|^2 - (y_i - x_i)(y_j - x_j)}{\| y-x \|^3} \mathbf{u}_i \otimes \mathbf{u}_j.
\end{align*}
The adjoint is obtained simply by commuting $\mathbf{u}_i$ and $\mathbf{u}_j$, and the composition of the adjoint with $D_y \pi_x |_{T_y \Sigma}$ is therefore
\begin{align*}
    & \big( D_y \pi_x |_{T_y \Sigma} \big)^* \big( D_y \pi_x |_{T_y \Sigma} \big) = \left( \sum_{i=1}^{n-1} \sum_{k=1}^n \frac{\delta_{ik} \| y-x \|^2 - (y_i - x_i)(y_k - x_k)}{\| y-x \|^3} \mathbf{u}_i \otimes \mathbf{u}_k \right) \\
    &\hs{4.5} \circ \left( \sum_{k=1}^n \sum_{j=1}^{n-1} \frac{\delta_{kj} \| y-x \|^2 + (y_k - x_k)(y_j - x_j)}{\| y-x \|^3} \mathbf{u}_k \otimes \mathbf{u}_j \right) \\
    &\hs{0.4} = \frac{1}{\| y-x \|^6} \sum_{i=1}^{n-1} \sum_{k=1}^n \sum_{j=1}^{n-1} \big( \delta_{ij} \| y-x \|^4 - \delta_{ik} (y_k - x_k)(y_j - x_j) \| y-x \|^2 \\
    &\hs{2.2} - \delta_{kj} (y_i - x_i)(y_k - x_k) \| y-x \|^2 + \delta_{ij} (y_i - x_i)(y_k - x_k)^2(y_j - x_j) \big) \mathbf{u}_i \otimes \mathbf{u}_j.
\end{align*}

Reordering the basis if necessary so that $y_1 - x_1 \neq 0$, we find that the eigenvectors of this operator are
\begin{align*}
    & (y_1 - x_1) \mathbf{u}_i - (y_i - x_i) \mathbf{u}_1, \quad i = 2, ..., n-1, \quad \text{with eigenvalue} \quad \frac{1}{\| y-x \|^2} \quad \text{and} \\
    & \sum_{i=1}^{n-1} (y_i - x_i) \mathbf{u}_i \quad \text{with eigenvalue} \quad \frac{\| y'-x' \|^2 \| y-x \|^2 - 2 \| y'-x' \|^2 \| y-x \|^2 + \| y-x \|^4}{\| y-x \|^6},
\end{align*}
where $z = (z',z_n)$. This last eigenvalue simplifies to
\begin{align*}
    &\frac{\| x'-y' \|^2 - 2 \| x'-y' \|^2 + \| x-y \|^2}{\| x-y \|^4} = \frac{\| x-y \|^2 - \| x'-y' \|^2}{\| x-y \|^4} \\
    &\hs{1.5} = \frac{(x_n - y_n)^2}{\| x-y \|^4} = \frac{\langle x-y, \mathbf{u}_n \rangle^2}{\| x-y \|^4} = \frac{\langle x-y, \nu(y) \rangle^2}{\| x-y \|^4},
\end{align*}
and taking the square root of the product of the eigenvalues yields the Jacobian:
\begin{equation*}
    |J\pi_x(y)| = \left( \prod_{i=1}^{n-2} \frac{1}{\| x-y \|^2} \right)^{\!1/2} \left( \frac{\langle x-y, \nu(y) \rangle^2}{\| x-y \|^4} \right)^{\!1/2} = \frac{|\langle x-y, \nu(y) \rangle|}{\| x-y \|^n}. \qedtag
\end{equation*}

\vs{0.15}

\define{Proof of the Theorem.} Let $x \in \Sigma$ be a point at which the OCC is satisfied. We employ the coarea formula by radially projecting $\Sigma$ onto $\partial B(x,1) \cong \bbs^{n-1}$ and applying Lemma \ref{lem:jacobian}:
\begin{align*}
    \int_\Sigma & \frac{\langle x-y, \nu(y) \rangle \langle y-x, \nu(x) \rangle}{\| x-y \|^{n+1}} \+ d\mathcal{H}^{n-1}(y) \\
    &\hs{1.5}= \int_\Sigma \frac{\langle y-x, \nu(x) \rangle}{\|x-y\|} |J\pi_x(y)| \sgn \, \langle x-y, \nu(y) \rangle \, d\mathcal{H}^{n-1}(y) \\
    &\hs{1.5}= \int_{\bbs^{n-1}} \int_{\Sigma \+ \cap \+ \pi_x^{-1}(\omega)} \left\langle \frac{y-x}{\|y-x\|}, \nu(x) \right\rangle \sgn \, \langle -\omega, \nu(y) \rangle \, d\mathcal{H}^0(y) \+ d\mathcal{H}^{n-1}(\omega) \\
    &\hs{1.5}= \int_{\bbs^{n-1}} \sum_{\substack{\lambda \+ > \+ 0 \\ x+\lambda\omega \+ \in \+ \Sigma}} \langle \omega, \nu(x) \rangle \sgn \, \langle -\omega, \nu(x + \lambda\omega) \rangle \, d\mathcal{H}^{n-1}(\omega) \\
    &\hs{1.5}= \frac{1}{2} \int_{\bbs^{n-1}} \sum_{y \in \Sigma \+ \cap \+ (L_{x,\omega} \setminus \{ x \})} \langle -\omega, \nu(x) \rangle \sgn \, \langle \omega, \nu(y) \rangle \, d\mathcal{H}^{n-1}(\omega),
  \end{align*}
where $L_{x,\omega}$ is the line through $x$ with direction vector $\pm\omega$. By the orientation cancellation condition (Equation \eqref{eq:occ}),
\begin{equation*}
    \sum_{y \+ \in \+ \Sigma \+ \cap \+ (L_{x,\omega} \setminus \{ x \})} \sgn \, \langle \omega, \nu(y) \rangle = \left( \sum_{y \+ \in \+ \Sigma \+ \cap \+ L_{x,\omega}} \sgn \, \langle \omega, \nu(y) \rangle \right) - \sgn \, \langle \omega, \nu(x) \rangle = \sgn \, \langle -\omega, \nu(x) \rangle,
\end{equation*}

so we conclude that
\begin{equation} \label{eq:apply-occ}
\begin{aligned}
    & \int_\Sigma \frac{\langle x-y, \nu(y) \rangle \langle y-x, \nu(x) \rangle}{\| x-y \|^{n+1}} \+ d\mathcal{H}^{n-1}(y) \\
    &\quad = \frac{1}{2} \int_{\bbs^{n-1}} \langle -\omega, \nu(x) \rangle \sgn \, \langle -\omega, \nu(x) \rangle \, d\mathcal{H}^{n-1}(\omega) \\
    &\quad = \frac{1}{2} \int_{\bbs^{n-1}} |\langle \omega, \nu(x) \rangle| \+ d\mathcal{H}^{n-1}(\omega).
\end{aligned}
\end{equation}
This last integral is invariant under rotations and may therefore be computed in graph coordinates on $\bbs^{n-1}$ with $\nu(x) = (0, ..., 0, 1)$, giving \vs{0.1}
\begin{equation*}
    \frac{1}{2} \int_{\bbs^{n-1}} |\langle \omega, \nu(x) \rangle| \+ d\mathcal{H}^{n-1}(\omega) = \frac{1}{2} \int_{\bbs^{n-1}} |\omega_n| \, d\mathcal{H}^{n-1}(\omega) = \frac{1}{2} \int_{B(0,1)} 2 \, d\mathcal{L}^{n-1}(x) = \alpha_{n-1}. \vs{0.1}
\end{equation*}
This combines with Equation \eqref{eq:apply-occ} to yield Equation \eqref{eq:main-1}, and since $\Sigma$ satisfies the OCC, this conclusion holds for $\mathcal{H}^{n-1}$-a.e.~$x \in \Sigma$. Equation \eqref{eq:main-2} then follows by integrating over $\Sigma$ with respect to $d\mathcal{H}^{n-1}(x)$. \hfill $\square$

Again, we precede the proof of the Corollary with a technical lemma to the effect that, when entering or leaving a set of finite perimeter, one typically must cross the reduced boundary. Denote by $E^{(t)}$ the set of points in $\R^n$ whose Lebesgue density with respect to $E$ is $t$ ($0 \leq t \leq 1$).

\begin{lem} \label{lem:finite-perimeter}
    If $E$ is a bounded set of finite perimeter, then $(\partial^*\-E,\nu_E)$ satisfies the OCC.
\end{lem}

\textit{Proof.} Modifying $E$ on an $\mathcal{L}^n$-null set\textemdash an operation that does not affect $\partial^*\-E$\textemdash we assume without loss of generality that $E^{(1)} \subseteq E$ and $E^{(0)} \subseteq \R^n \setminus E$. Let $\chi_E^*$ be the precise representative of $\chi_E$ and, for $x \in \R^n$ and $\omega \in \bbs^{n-1}$, define $(\chi_E^*)_x^\omega \!: \R \to \R$ by
\begin{equation} \label{eq:precise-rep}
    (\chi_E^*)_x^\omega(t) := \chi_E^*(x + t\omega).
\end{equation}
Our claim will follow from a general result, adapted here from \cite{ambrosio2000functions} Theorem 3.108:
\begin{quote}
    \textbf{Theorem.} \textit{For $\mathcal{H}^{n-1}$-a.e.~$x \in \R^n$, the following statements hold for $\mathcal{H}^{n-1}$-a.e.~$\omega \in \bbs^{n-1}$: \vs{0.1}
    \begin{enumerate}[label={\normalfont \arabic*}.,noitemsep]
        \item $L_{x,\omega} \subseteq E^{(0)} \cup \partial^*\-E \cup E^{(1)}$.
        \item $(\chi_E^*)_x^\omega$ has bounded variation, $(\chi_E^*)_x^\omega(t) = \chi_E(x + t\omega)$ for $\mathcal{L}^1$-a.e.~$t \in \R$, and the set of discontinuities of $(\chi_E^*)_x^\omega$ is given by
        \begin{equation*}
            (J_E)_x^\omega := \{ t \in \R \!: x + t\omega \in \partial^*\-E \}.
        \end{equation*}
        \item $\langle \omega, \nu_E(y) \rangle \neq 0$ for all $y \in \partial^*\-E \cap L_{x,\omega}$.
        \item For all $x + t\omega \in \partial^*\-E$,
        \begin{align*}
            & \lim_{s \+ \uparrow \+ t} \ (\chi_E^*)_x^\omega(s) = \left\{ \begin{array}{cl}
                1 & \text{if } \langle \omega, \nu_E(x + t\omega) \rangle > 0 \\
                0 & \text{if } \langle \omega, \nu_E(x + t\omega) \rangle < 0
            \end{array}\right. \quad \text{and} \\
            & \lim_{s \+ \downarrow \+ t} \ (\chi_E^*)_x^\omega(s) = \left\{ \begin{array}{cl}
                0 & \text{if } \langle \omega, \nu_E(x + t\omega) \rangle > 0 \\
                1 & \text{if } \langle \omega, \nu_E(x + t\omega) \rangle < 0.
            \end{array}\right.
        \end{align*}
    \end{enumerate}}
    \vs{0.1}
\end{quote}
Claims 1 and 3 are included to make sense of Claims 2 and 4, respectively. By the structure of bounded sets of finite perimeter on the real line, Claim 2 also implies that $(\chi_E^*)_x^\omega$ is equal to the indicator function of a finite collection of positively separated bounded intervals, except at the endpoints of these intervals (the points of $(J_E)_x^\omega$), where $(\chi_E^*)_x^\omega$ takes the value $\tfrac{1}{2}$. (See \cite{maggi2012sets} Proposition 12.13 and \cite{ambrosio2000functions} Theorem 3.28.)

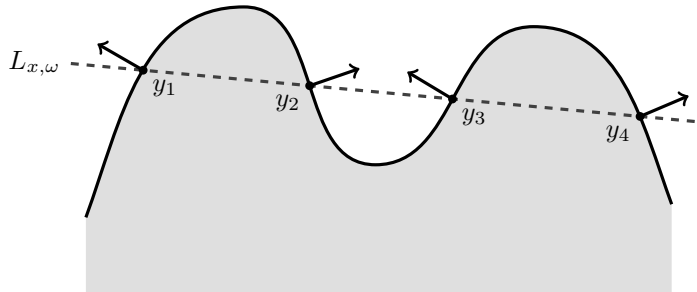
\begin{figure}
    \centering
    \begin{tikzpicture}[scale=0.875]
    
    \fill[verylightgray] (-4.5,-2.25) -- (-4.5,-1.1) to [out=70,in=180] (-2.1,2.1) to [out=0,in=180] (-0.1,-0.3) to [out=0,in=180] (2.3,1.8) to [out=0,in=110] (4.4,-0.9) -- (4.4,-2.25) -- (-4.5,-2.25);
    \draw[very thick] (-4.5,-1.1) to [out=70,in=180] (-2.1,2.1) to [out=0,in=180] (-0.1,-0.3) to [out=0,in=180] (2.3,1.8) to [out=0,in=110] (4.4,-0.9);
    
    \draw[very thick,dashed,darkgray] (4.75,0.35814) -- (-4.75,1.24186);

    \normvec{-3.63}{1.138}{0.8}{150};
    \normvec{-1.1}{0.902}{0.8}{19};
    \normvec{1.07}{0.7}{0.8}{149};
    \normvec{3.92}{0.435}{0.8}{23};
    
    \fill (-3.63,1.138) circle (0.67mm);
    \fill (-1.1,0.902) circle (0.67mm);
    \fill (1.07,0.70) circle (0.67mm);
    \fill (3.92,0.435) circle (0.67mm);

    \node[left] at (-4.75,1.24186) {$L_{x,\omega}$};
    
    \node[below right] at (-3.63,1.138) {$y_1$};
    \node[below left] at (-1.1,0.902) {$y_2$};
    \node[below right] at (1.07, 0.70) {$y_3$};
    \node[below left] at (3.92,0.435) {$y_4$};
    
    \end{tikzpicture}

    \vs{0.3}
    
    \caption{For almost all $x \in \partial^*\-E$ and $\omega \in \bbs^{n-1}$, the arrangement of the normal vectors to $\partial^*\-E$ along $L_{x,\omega}$ is the ``obvious" one depicted here. This is the thrust of Lemma \ref{lem:finite-perimeter} and the reason Equation \eqref{eq:angles} implies convexity in the proof of the Corollary.}
    \label{fig:in-n-out}
\end{figure}

Let $x$ and $\omega$ be as above and enumerate the points of $(J_E)_x^\omega$ by $y_i = x + t_i \omega$, $t_1 < \cdots < t_k$. The preceding remark on the structure of $(\chi_E^*)_x^\omega$ implies that $k$ is even, so Equation \eqref{eq:occ} defining the OCC will be satisfied if $\langle \omega, \nu_E(y_i) \rangle$ and $\langle \omega, \nu_E(y_{i+1}) \rangle$ have opposite signs for $i = 1, ..., k-1$, as in Figure \ref{fig:in-n-out}. If $\langle \omega, \nu_E(y_i) \rangle > 0$, then, by Claim 4,
\begin{equation*}
    \lim_{s \+ \downarrow \+ t_i} \ (\chi_E^*)_x^\omega(s) = 0.
\end{equation*}
By Claim 2 and the structure of $(\chi_E^*)_x^\omega$, we must have $\chi_E^*(x + t\omega) = 0$ for all $t > t_i$ sufficiently close to $t_i$. Since $(\chi_E^*)_x^\omega$ is continuous on $(t_i,t_{i+1})$ and takes a discrete set of values, it follows that $\chi_E^*(x + t\omega) = 0$ for all $t_i < t < t_{i+1}$, whence
\begin{equation*}
    \lim_{s \+ \uparrow \+ t_{i+1}} \ \chi_E^*(x + t\omega) = \lim_{s \+ \uparrow \+ t_{i+1}} \ (\chi_E^*)_x^\omega(s) = 0.
\end{equation*}

Another application of Claim 4 gives $\langle \omega, \nu_E(y_{i+1}) \rangle < 0$, as desired. By an identical argument, we have $\langle \omega, \nu_E(y_{i+1}) \rangle > 0$ whenever $\langle \omega, \nu_E(y_i) \rangle < 0$, so the signs of the inner products against $\omega$ alternate, as we sought to show. \hfill $\square$ \vs{0.15}

Theorem 3.108 of \cite{ambrosio2000functions} (which the authors attribute to Vol'pert \cite{vol1967spaces}) actually implies a little more: for \textit{every} $\omega \in \bbs^{n-1}$, the set of $x \in \partial^*\-E$ such that Equation \eqref{eq:occ} fails (with $(\Sigma,\nu) = (\partial^*\-E,\nu_E)$) is $\mathcal{H}^{n-1}$-null. However, this leaves open the question whether it is also true that, for every $x \in \partial^*\-E$, the set of $\omega \in \bbs^{n-1}$ such that Equation \eqref{eq:occ} fails is $\mathcal{H}^{n-1}$-null. \vs{0.15}

\define{Proof of the Corollary.} By Lemma \ref{lem:finite-perimeter}, $(\Sigma,\nu) = (\partial^*\-E,\nu_E)$ satisfies Equation \eqref{eq:main-2}, and applying the triangle inequality for integrals gives \eqref{eq:convex-inequality}. Equality holds if and only if the integrand in Equation \eqref{eq:main-2} is almost everywhere nonnegative, i.e., if and only if
\begin{equation} \label{eq:angles}
    \langle x-y, \nu_E(y) \rangle \langle y-x, \nu_E(x) \rangle \geq 0 \qquad \text{for } \mathcal{H}^{n-1}\text{-a.e.~} x,y \in \partial^*\-E.
\end{equation}
Recall from the proof of Lemma \ref{lem:finite-perimeter} that sets of finite perimeter satisfy the OCC in a strong way; namely, a typical line in $\R^n$ parametrized by $t \mapsto x + t\omega$ intersects $\partial^*\-E$ at points $y_i = x + t_i \omega$, $t_1 < \cdots < t_k$, such that $\langle \omega, \nu_E(y_i) \rangle$ and $\langle \omega, \nu_E(y_{i+1}) \rangle$ have opposite signs for all $i$. In particular, since $\chi_E(x + t\omega) = 0$ for $\mathcal{L}^1$-a.e.~$t \in (-\infty,t_1)$, Claim 4 above entails that $\langle \omega, \nu_E(y_1) \rangle < 0$, from which it follows that $\langle \omega, \nu_E(y_i) \rangle < 0$ for $i$ odd and $\langle \omega, \nu_E(y_i) \rangle > 0$ for $i$ even. With $y_i - y_1 = \| y_i - y_1 \| \+ \omega$ for any $i > 1$, we get
\begin{equation*}
    \langle y_i - y_1, \nu_E(y_1) \rangle < 0
\end{equation*}
and, under the hypothesis of Equation \eqref{eq:angles},
\begin{equation*}
    \langle y_1 - y_i, \nu_E(y_i) \rangle < 0, \quad \text{i.e.,} \quad \langle \omega, \nu_E(y_i) \rangle > 0.
\end{equation*}
But then $i > 1$ is even yet arbitrary, so it must be that $k = i = 2$. That is, for a.e.~$x \in \partial^*\-E$ and $\omega \in \bbs^{n-1}$, the line $L_{x,\omega}$ intersects $\partial^*\-E$ in exactly two points, from which it follows that either $E$ or $\R^n \setminus E$ is $\mathcal{L}^n$-equivalent to a convex set. Because $E$ is bounded, it must be the former.

Conversely, if $E$ is $\mathcal{L}^n$-equivalent to a convex set, then we modify it without loss of generality on a null set so that it is (truly) convex. This ensures that $\partial^*\-E$ is a dense subset of $\partial E$, so we may write
\begin{equation*}
    \overline{E} = \bigcap_{x \in \partial^*\-E} \{ y \in \R^n \!: \langle y-x, \nu_E(x) \rangle \leq 0 \}.
\end{equation*}
Consequently, for all $x,y \in \partial^*\-E \subseteq \overline{E}$, we have both $\langle x-y, \nu_E(y) \rangle \leq 0$ and $\langle y-x, \nu_E(x) \rangle \leq 0$, from which it follows that the integrand Equation \eqref{eq:main-2} is nonnegative and that \eqref{eq:convex-inequality} holds with equality. \hfill $\square$

\phantomsection
\section*{Acknowledgement}
My gratitude goes to Stefan Steinerberger for his ideas and suggestions throughout the drafting of this article.

\bibliographystyle{plain}
\bibliography{references}

\end{document}